\documentclass{article}
\usepackage{spconf,amsmath,graphicx}

\usepackage{amssymb}  
\usepackage[utf8]{inputenc}
\usepackage{color}
\usepackage{tikz}
\usetikzlibrary{shapes,arrows.meta}
\usepackage{siunitx}

\usepackage[font=footnotesize,labelfont=bf]{caption}

\usepackage{algorithm}
\usepackage{algpseudocode}
\usepackage{cite}
\usepackage{bm,mathabx,fixmath,}

\algrenewcommand\algorithmicrequire{\textbf{Input:}}

\DeclareMathOperator*{\argmin}{arg\,min}

\DeclareMathOperator{\prox}{prox}

\newtheorem{remark}{Remark}

\newtheorem{theorem}{Theorem}
\newtheorem{assumption}{Assumption}

\newcommand{\norm}[1]{\left\lVert#1\right\rVert}
\newcommand{\x}{\mathbold{x}}
\newcommand{\y}{\mathbold{y}}
\newcommand{\w}{\mathbold{w}}
\newcommand{\z}{\mathbold{z}}
\newcommand{\p}{\mathbold{p}}
\newcommand{\R}{\mathbb{R}}
\newcommand{\Q}{\mathbold{Q}}

\newcommand{\env}{\operatorname{M}}
\newcommand{\afbe}{\tilde{\operatorname{M}}}

\newcommand{\svmin}[1]{\sigma_m\left(#1\right)}
\newcommand{\svmax}[1]{\sigma_M\left(#1\right)}

\def\endstatement{\hfill$\square$}

\newcommand{\A}{\mathbf{A}}
\newcommand{\bv}{\mathbf{b}}

\title{Prediction-Correction for Nonsmooth Time-Varying Optimization via Forward-Backward Envelopes}
\name{Nicola Bastianello$^*$, Andrea Simonetto$^\S$, Ruggero Carli$^*$\vspace*{-3mm}}
\address{$^*$Department of Information Engineering, University of Padova 
  \\ $^\S$IBM Research Ireland \vspace*{-3mm}}

\begin{document}
\ninept
\maketitle
%

\begin{abstract}
We present an algorithm for minimizing the sum of a strongly convex time-varying function with a time-invariant, convex, and nonsmooth function.
The proposed algorithm employs the prediction-correction scheme alongside the forward-backward envelope, and we are able to prove the convergence of the solutions to a neighborhood of the optimizer that depends on the sampling time.
Numerical simulations for a time-varying regression problem with elastic net regularization highlight the effectiveness of the algorithm.
\end{abstract}


\section{Introduction}\label{sec:intro}
In this work, we are interested in the solution of time-varying optimization problems in the form
\begin{equation}\label{eq:tv-problem}
	\x^*(t) = \argmin_{\x\in \R^n} \{ f(\x;t) + g(\x) \}
\end{equation}
where $f:\R^n \times \R_+ \to \R$ is smooth and strongly convex, and $g:\R^n \to \R$ is proper, closed and convex, but possibly non-differentiable. Since the solution $\x^*(t)$ -- the trajectory -- changes over time, the objective is to track it up to a bounded error ball.

In particular, we are interested in solving problem \eqref{eq:tv-problem} in a discrete-time framework, in order to directly implement the solution on digital hardware. Hence we discretize the problem with a sampling period $T_{\mathrm{s}}:= t_{k+1} - t_{k}$, which yields the sequence of time-invariant problems
\begin{equation}\label{eq:ti-problem}
	\x^*(t_{k+1}) = \argmin_{\x\in \R^n} \{ f(\x;t_{k+1}) + g(\x) \}, \quad k \in \mathbb{N}.
\end{equation}

The smaller the sampling time is, the higher the accuracy of the trajectory composed of the solutions to~\eqref{eq:ti-problem} will be. However, we need to account for the time required to solve the problems which might exceed some values of $T_\mathrm{s}$, and therefore there is a trade-off between precision and practical implementation constraints.

There are many applications in which problems in the form~\eqref{eq:tv-problem} arise. For instance, in signal processing the reconstruction of time-varying signals on the basis of (noisy) observations gathered online can be cast as a sequence of optimization problems \cite{asif2010dynamic,asif2014sparse,Vaswani2015,Yang2015,charles2016dynamic,Sopasakis2016}. In control, the \textit{model predictive control} (MPC) requires that we solve an optimization problem which varies over time \cite{jerez2014embedded,hours2016parametric,gutjahr2017lateral} in order to design a control action. In robotics, path tracking and leader following problems can be cast in the framework of \eqref{eq:tv-problem}, see for example \cite{verscheure2009time,ardeshiri2011convex,dixit2018online}.

In this paper, we are interested in the solution of~\eqref{eq:ti-problem} using a \emph{prediction-correction} scheme. Time-varying optimization algorithms based on the prediction-correction scheme have been proposed for both the discrete-time framework that we employ \cite{simonetto2016class,simonetto2017prediction,simonetto2018tac} and in a continuous-time setup \cite{fazlyab2017prediction,Rahili2015,fazlyab2016self}.

These works, however, are designed to solve smooth optimization problems only; here, our aim is to tackle non-smooth optimization problems by employing the recent results on \textit{envelope functions}, and in particular the \textit{forward-backward envelope} (FBE) first introduced in \cite{patrinos2013proximal}, in conjunction with the prediction-correction scheme.

The main contributions of this paper are: \emph{(i)} a prediction-correction algorithm to solve the time-varying optimization problem~\eqref{eq:tv-problem} by using the envelope functions in both the prediction and correction step; \emph{(ii)} a detailed convergence and convergence rate analysis of the above that show global convergence to an error bound of $O(T_{\mathrm{s}})$ and local convergence to an error bound of $O(T_{\mathrm{s}}^2)$, under additional assumptions.

\begin{remark}
The forward-backward envelope is a powerful tool that has recently gained momentum, especially in the context of solving certain classes on non-convex optimization problems. In this sense, this paper can be seen as a first step towards a more general theory of time-varying optimization algorithms. We remark also that the FBE has been advocated as a way to derive Newton-like methods for $\ell_1$ regularized problems, showing improved results in comparison to more traditional approaches, such as FISTA~\cite{Sopasakis2016}, at lower computational cost. In Sec.~\ref{sec:background}, we will report some results about the FBE, and we refer the reader to \cite{stella2017forward,themelis2018forward,giselsson2018envelope} for an in-depth treatment of the subject. 
\end{remark}

\begin{remark}
For the relationship of the FBE with the \textit{forward-backward splitting} (also known as \textit{proximal gradient method}) \cite{parikh_proximal_2014,combettes_proximal_2011}, see \textit{e.g.} \cite{stella2017forward}.
\end{remark}

\paragraph*{Organization} The paper is organized as follows. Sec.~\ref{sec:background} introduces the prediction-correction scheme, the forward-backward envelope, and then the proposed algorithm. Sec.~\ref{sec:convergence} presents the convergence results for the algorithm and a sketch of the proof. Sec.~\ref{sec:simulation} describes the results of the numerical simulations and Sec.~\ref{sec:conclusions} some concluding remarks.

\paragraph*{Basic definitions} We say that a function $\varphi:\R^n \to \R$ is $m$-\textit{strongly convex} for a constant $m \in \R_+$ iff $\varphi(\x) - \frac{m}{2}\|\x\|^2$ is convex. The function $\varphi$ is said to be $L$-\textit{smooth} if its gradient is $L$-Lipschitz continuous, or equivalently $\varphi(\x) - \frac{L}{2}\|\x\|^2$ is concave. We denote the class of $m$-strongly convex and $L$-smooth functions with $\mathcal{S}_{m,L}(\R^n)$.

A function is said to be \textit{closed} if for any $a \in \R$ the set $\{ \x \in \operatorname{dom}(f)\ |\ \varphi(\x) \leq a \}$ is closed. A function is said to be \textit{proper} if it does not attain $-\infty$. We denote the class of closed, convex and proper functions with $\Gamma_0(\R^n)$.

Given $\varphi \in \Gamma_0(\R^n)$ we define its \textit{subdifferential} as the set-valued operator $\partial \varphi: \R^n \rightrightarrows \R^n$ such that
$$
	\x \mapsto \left\{ \z \in \R^n\ |\ \forall \y \in \R^n:\ \langle \y-\x, \z \rangle + \varphi(\x) \leq \varphi(\y) \right\}.
$$

\section{Prediction-Correction with Envelopes}\label{sec:background}
In this section, we introduce the prediction-correction scheme for time-varying optimization alongside with the forward-backward envelope function. In the remainder of this paper we make use of the following assumptions.

\begin{assumption}\label{as:first}
The function $f:\R^n \times \R_+ \to \R^n$ belongs to $\mathcal{S}_{m,L}(\R^n)$ uniformly in time.  The function $g:\R^n \to \R$ belongs to $\Gamma_0(\R^n)$ and is in general nonsmooth. 
\end{assumption}

\begin{assumption}\label{as:first-bis}
The function $f$ has bounded time derivative of its gradient derivative as: 
$\norm{\nabla_{t\x} f(\x;t)} \leq C_0$.
\end{assumption}

\begin{assumption}\label{as:second}
The function $f$ is at least three time differentiable and has bounded derivatives w.r.t. $\x \in \R^n$ and $t \in \R_+$ as: 
\begin{align*}
	\norm{\nabla_{\x\x\x} f(\x;t)} & \leq C_1, \quad \norm{\nabla_{\x t \x} f(\x;t)} \leq C_2, \\ &\norm{\nabla_{tt\x} f(\x;t)} \leq C_3.
\end{align*}
\end{assumption}

\smallskip

In the analysis of time-varying problems, Assumption~\ref{as:first} is common, see \textit{e.g.} \cite{simonetto2016class,popkov2005gradient,dontchev2013euler}. This assumption ensures by strong convexity that the solution to the problem is unique at each time, and that the gradient of $f$ is Lipschitz continuous. Moreover, Assumption~\ref{as:first-bis} guarantees that the gradient of $f$ has a variability over time that is bounded, thus enabling the computation of reliable predictions. Assumption~\ref{as:second} imposes instead boundedness of the tensor $\nabla_{\x\x\x}f(\x;t)$, which is typical when analyzing the convergence of second-order algorithms. Moreover, it bounds the variability of the Hessian of $f$ over time, which makes it possible to carry out even more precise predictions of the optimal trajectory.

\subsection{Prediction-correction}

Prediction-correction algorithms have appeared as a computational-light way to solve time-varying optimization problems. The main idea is to compute approximate optimizers for the sequence of time-invariant problems~\eqref{eq:ti-problem}, such that eventually one converges on the time-varying optimizer trajectory $\x^*(t)$. More formally, let $\x_k$ be the approximate optimizer for \eqref{eq:ti-problem} at $k$. Then we want to design methods to determine the sequence $\{\x_k\}_{k\in \mathbb{N}}$ such that $\|\x_k - \x^*(t_k)\|$ goes eventually to a bounded error term. 

Prediction-correction algorithms determine each $\x_{k+1}$ by first predicting (at $t_{k}$) how the optimizer will change in time, and then by correcting (at $t_{k+1}$) based on the new acquired sampled cost function. Both prediction and correction are here based on a few descent iterations on the envelope functions. The more iterations one performs, the smaller the asymptotical tracking error, however the greater the computational time is. 


The \textit{prediction} step has the aim of computing an approximation of the optimal solution at time $t_{k+1}$, $\x^*(t_{k+1}) =: \x_{k+1}^*$, by using only the information available at time $t_k$, that is $f(\x;t_k)$ and the previous solution $\x_k$ computed by the algorithm. Once the new cost function $f(\x;t_{k+1})$ is observed at time $t_{k+1}$, we perform the \textit{correction} step, that is we solve problem~\eqref{eq:ti-problem} approximately, using as initial condition the prediction computed at time $t_k$.

In order to use the forward-backward envelope framework, it is useful to reformulate the minimization problem~\eqref{eq:ti-problem} as the following generalized equation 
\begin{equation}\label{eq:gen-equation}
	\nabla_{\x} f(\x_{k+1};t_{k+1}) + \partial g(\x_{k+1}) \ni 0.
\end{equation}

During the prediction step at time $t_{k}$, we cannot solve~\eqref{eq:gen-equation} to predict how the optimizer will change at $t_{k+1}$; instead, we make use of the available information at time $t_k$ to approximate $\nabla_{\x} f(\x;t_{k+1})$ with the following Taylor expansion
\begin{align}
\begin{split}
	\nabla h_k(\x) &= \nabla_{\x} f(\x_k;t_k) + \\ & + \nabla_{\x\x} f(\x_k;t_k) (\x - \x_k) + T_{\mathrm{s}} \nabla_{t\x} f(\x_k;t_k).
\end{split}
\end{align}
Therefore during the prediction step we want to solve the approximated generalized equation
\begin{equation}\label{eq:approx-gen-equation}
	\nabla h_k(\x_{k+1|k}) + \partial g(\x_{k+1|k}) \ni 0
\end{equation}
derived from~\eqref{eq:gen-equation} substituting $h_k(\x)$ to the (as yet unknown) $f(\x;t_{k+1})$; notice that $\x_{k+1|k}$ will denote the computed prediction.

During the correction step at time $t_{k+1}$, we can now solve (approximately)~\eqref{eq:gen-equation}, which is what we will do.  

\begin{remark}
From Assumption~\ref{as:first} follows that $h_k \in \mathcal{S}_{m,L}(\R^n)$, and by definition we can write it explicitly as
\begin{align*}
	h_k(\x) &= \frac{1}{2} \x^\top \nabla_{\x\x} f(\x_k;t_k) \x + \\ +& \Big( \nabla_{\x} f(\x_k;t_k) - \nabla_{\x\x} f(\x_k;t_k) \x_k + T_{\mathrm{s}} \nabla_{t\x} f(\x_k;t_k) \Big)^\top \x.
\end{align*}
\end{remark}

\subsection{Forward-backward envelope}
Notice that both the prediction and correction problems, \eqref{eq:approx-gen-equation} and \eqref{eq:gen-equation}, are of the form
\begin{equation}\label{eq:gen-problem}
	\nabla \varphi(\x^*) + \partial g(\x^*) \ni 0
\end{equation}
with $\varphi$ that is $m$-strongly convex and $L$-smooth. Therefore we can apply the recently proposed \textit{forward-backward envelope} (FBE) to solve them.

The FBE for a problem \eqref{eq:gen-problem} is defined as
\begin{equation}\label{eq:fbe}
	\env(\x) = \min_\y \Bigg\{ \varphi(\x) + \langle \nabla \varphi(\x), \y - \x \rangle + g(\y) + \frac{\norm{\y - \x}^2}{2\gamma} \Bigg\}
\end{equation}
where $\gamma \in (0, 1/L)$.

Under Assumption~\ref{as:first} it holds that 
$$
\argmin (\varphi + g)(\x) = \argmin \env(\x),
$$ 
and therefore minimizing the FBE is equivalent to solving problem \eqref{eq:gen-problem}. Moreover, the envelope is continuously differentiable on $\R^n$ and twice continuously differentiable at the unique solution $\x^*$, with positive definite Hessian.

In general, the FBE is however nonconvex, and hence in order to minimize it a \textit{quasi-Newton} scheme with line search has been proposed in \cite{stella2017forward}, that estimates the Hessian of the FBE using the BFGS method. In our framework, it is possible to prove that the quasi-Newton method applied to the FBE has global linear convergence, that is
$$
	\norm{\x^{l+1} - \x^*} \leq \zeta \norm{\x^l - \x^*}, \quad l \in \mathbb{N}
$$
with
\begin{equation}\label{eq:zeta}
	\zeta = \sqrt{\max\left\{ \frac{1}{2}, 1 - \frac{m}{4} \min\left\{ \gamma, \frac{1}{4L} \right\} \right\}} < 1,
\end{equation}
a result that will be instrumental in proving convergence of our prediction-correction algorithm.

\begin{remark}\label{rem:convexity-fbe}
The recent work \cite{giselsson2018envelope} proved that if $\varphi$ is convex quadratic, then the FBE is strongly convex and smooth, and notice that this is exactly the case of $h_k$ in the prediction step. Therefore we can minimize the FBE at the prediction step using a Newton method with BFGS scheme, without the need for the line search that requires a larger number of iterations. The numerical results presented in Sec.~\ref{sec:simulation} exploit this.
\end{remark}

\begin{remark}
An alternative minimization strategy for the FBE is proposed in \cite{themelis2018forward}.
\end{remark}

\subsection{Proposed algorithm}
The previous section introduced the forward-backward envelope, that is suited to solving the prediction and correction problems. However, the convergence of the quasi-Newton method is guaranteed only asymptotically. For practical reasons, namely the finite length of each sampling period, we choose to perform only a fixed number of iterations of the solution algorithm: $P$ for the prediction step, $C$ for the correction.

We are now ready to describe the proposed prediction-correction algorithm with the FBE, which is reported in Algorithm~\ref{alg:envelope-prediction-correction}: at every time $t_k$, we perform $P$ steps of the quasi-Newton method for the FBE, 
\begin{equation*}
	\afbe(\x)\! =\! \min_\y \Bigg\{ h_k(\x) +\langle \nabla h_k(\x), \y - \x \rangle + g(\y) + \frac{\norm{\y\! - \!\x}^2}{2\gamma} \Bigg\},
\end{equation*}
constructed for the prediction problem~\eqref{eq:approx-gen-equation} [cf. line~\ref{line1}]; this yields an approximate predictor $\tilde{\x}_{k+1|k}$. 

At time $t_{k+1}$, we observe the new cost function $f(\cdot;t_{k+1})$ [cf. line~\ref{line2}], and we perform $C$ steps of the quasi-Newton method for the FBE, 
\begin{multline*}
	\env(\x) = \min_\y \Bigg\{ f(\x;t_{k+1}) +\langle \nabla_{\x} f(\x; t_{k+1}), \y - \x \rangle +\\ + g(\y) + \frac{\norm{\y - \x}^2}{2\gamma}\Bigg\},
\end{multline*}
constructed for the correction problem~\ [cf. line~\ref{line3}]; this yields the approximate optimizer $\x_{k+1}$.  

\begin{algorithm}
\caption{Prediction-correction algorithm with the FBE.}
\label{alg:envelope-prediction-correction}
\begin{algorithmic}[1]
	\Require $\x_0$, parameter $\gamma$, horizons $P$ and $C$.
	\For{$k=0,1,\ldots$}
		\State // time $t_k$
		\State perform $P$ steps of the quasi-Newton method for the FBE with initial condition $\x_k$ \label{line1}
		\State set $\tilde{\x}_{k+1|k}$ equal to the last iterate produced by the quasi-Newton
		\State // time $t_{k+1}$
		\State observe the cost function $f(\cdot;t_{k+1})$ \label{line2}
		\State perform $C$ steps of the quasi-Newton method for the FBE with initial condition the prediction $\tilde{\x}_{k+1|k}$ \label{line3}
		\State set $\x_{k+1}$ equal to the last iterate produced by the quasi-Newton
	\EndFor
\end{algorithmic}
\end{algorithm}

\begin{remark}
In general we could use two different $\gamma$ parameters for the prediction and correction steps, but for simplicity we use a single one.
\end{remark}

\section{Convergence Analysis}\label{sec:convergence}
In this section, we prove that the sequence $\{\x_k\}_{k\in\mathbb{N}}$ generated by Algorithm~\ref{alg:envelope-prediction-correction} converges to a neighborhood of the optimal trajectory, which is characterized in terms of the sampling period $T_{\mathrm{s}}$. We divide the result in two theorems. The first is a global convergence result with standard assumptions; the second is a local enhanced convergence result with additional assumptions. Such results are typical in prediction-correction time-varying optimization and they extend the ones in~\cite{simonetto2017prediction} for non-smooth cost functions.  

\begin{theorem}\label{th:linear-convergence}
Let Assumptions~\ref{as:first}-\ref{as:first-bis} hold, and choose the parameters $P$ and $C$ in such a way that
$$
	\zeta^C \left[ \zeta^P + (\zeta^P + 1) \frac{2L}{m} \frac{1-\gamma m}{1 - \gamma L} \right] < 1.
$$
Then the trajectory $\{\x_k\}_{k \in \mathbb{N}}$ generated by Algorithm~\ref{alg:envelope-prediction-correction} converges to a neighborhood of the optimal trajectory $\{\x_k^*\}_{k \in \mathbb{N}}$ as
$$
	\limsup_{k \to \infty} \norm{\x_k - \x_k^*} = O(\zeta^C T_{\mathrm{s}}).
$$
\endstatement
\end{theorem}

\smallskip

\begin{theorem}\label{th:quadratic-convergence}
Let Assumptions~\ref{as:first},~\ref{as:first-bis},~and~\ref{as:second} hold, and choose the parameters $P$ and $C$, and $\tau \in (0,1)$ in such a way that $\zeta^{P+C} < \tau$.

\noindent Then there exist an upper bound for the sampling time $\bar{T}_{\mathrm{s}}$ and a convergence region $\bar{R}$ such that if $T_{\mathrm{s}} < \bar{T}_{\mathrm{s}}$ and $\norm{\x_0 - \x_0^*} < \bar{R}$, then
$$
	\limsup_{k \to \infty} \norm{\x_k - \x_k^*} = O(\zeta^C T_{\mathrm{s}}^2) + O(\zeta^{P+C} T_{\mathrm{s}}).
$$
In particular, the bound for the sampling time and the convergence region are characterized by
\begin{align*}
	\bar{T}_{\mathrm{s}} = \frac{\tau - \zeta^{P+C}}{\zeta^C (\zeta^P + 1)} \frac{1}{\kappa  ( \kappa C_0C_1 + C_2)}
\end{align*}
$$
	\bar{R} = \frac{2}{C_1} ( \kappa C_0C_1 + C_2 ) \left( \bar{T}_{\mathrm{s}} - \frac{m}{\zeta^C} T_{\mathrm{s}} \right).
$$
with $\kappa = (1-\gamma m)/ [m(1 - \gamma L)]$.
\endstatement
\end{theorem}

The two theorems guarantee that, under suitable regularity conditions of the problem in hand, the trajectory generated by Algorithm~\ref{alg:envelope-prediction-correction} converges asymptotically to a neighborhood of the optimal trajectory. Moreover, the size of this neighborhood depends on $T_\mathrm{s}$ for Theorem~\ref{th:linear-convergence} and on $T_\mathrm{s}^2$ for Theorem~\ref{th:quadratic-convergence}, in accordance with the fact that the smaller the sampling time is, the better the sequence of problems \eqref{eq:ti-problem} approximates the original problem \eqref{eq:tv-problem}.

The neighborhoods depend also on the convergence rate $\zeta$, which in turn depends on the convexity and smoothness moduli of the function $f$; thus the structure of the problem influences the accuracy of the proposed algorithm.


The proof of both results can be found in the Appendix, along with the exact expression for the asymptotic error. Here we mention only some facts. The idea behind the proof is to compute an upper bound to the error $\norm{\x_k - \x_k^*}$, and to do so we need to account for two sources of error: the approximation error introduced during the prediction step, and the early termination error due to the finite number of minimization steps in both prediction and correction. The approximation error depends (among other things) on how fast the cost function is changing and a bound on such error can be derived based on implicit function mapping theorems. In particular, we make use of Dini's theorem, see \textit{e.g.}, \cite[Th.~1B.1]{dontchev2014implicit}, and the algebraic properties of the envelope function. The early termination errors are instead bounded based solely on the properties of the envelope function. 

Once the bound for the errors is derived, we combine them and provide a bound for the error $\norm{\x_k - \x_k^*}$ based on the parameters of Algorithm~\ref{alg:envelope-prediction-correction} (\textit{i.e.}, the step-size $\gamma$, and the horizons $P$ and $C$). Then, we choose such parameters in order to guarantee a finite error. The error bound available is in general not tight, and therefore it might be possible to relax the conditions on the parameters while still ensuring the convergence; which we will explore in future research. 

\section{Simulations}\label{sec:simulation}

Inspired by~\cite{Sopasakis2016}, and only as a proof of concept of our algorithm, we consider a regression problem, where we are interested in reconstructing a sparse time-varying signal $\y_k$ from the noisy measurements $\bv_k = \A \y_k + \mathbf{e}_k$ where the matrix $\A \in \R^{m \times n}$, the measurement vector $\bv_k \in \R^m$ with $m < n$ -- in particular $m=25$ and $n=50$ -- and the components of the error vector $\mathbf{e}_k$ are drawn from the normal distribution $\mathcal{N}(0, 10^{-3})$. We apply an \textit{elastic net} to solve the problem
, \textit{i.e.}, we define $g(\x) = \alpha \norm{\x}_1$ and $f(\x;t_k) = (1/2) \norm{\A \x - \bv_k}^2_2 + (1-\alpha) \norm{\x}^2_2 /2$ for $\alpha \in [0,1]$, and we formulate the sequence of time-invariant problems
\begin{equation}
\x^*(t_{k}) = \argmin_{\x\in \R^n} \Bigg\{ \frac{1}{2} \norm{\A \x - \bv_k}^2_2 + \frac{(1-\alpha)}{2}\norm{\x}^2_2 + \alpha \norm{\x}_1 \Bigg\}.
\end{equation}

Each component of the signal to be reconstructed is either of the form $y_k^i = c \sin(\omega t_k + \phi^i)$ where $c, \phi^i$ are random, or it is $0$; $\omega$ is set as $1/20$, so that we do half-a-turn every minute. The number and index of the ``active'' components is fixed at $6$. 

\begin{figure}[htb]

\begin{minipage}[b]{1\linewidth}
  \centering
  \centerline{\includegraphics[width=\textwidth]{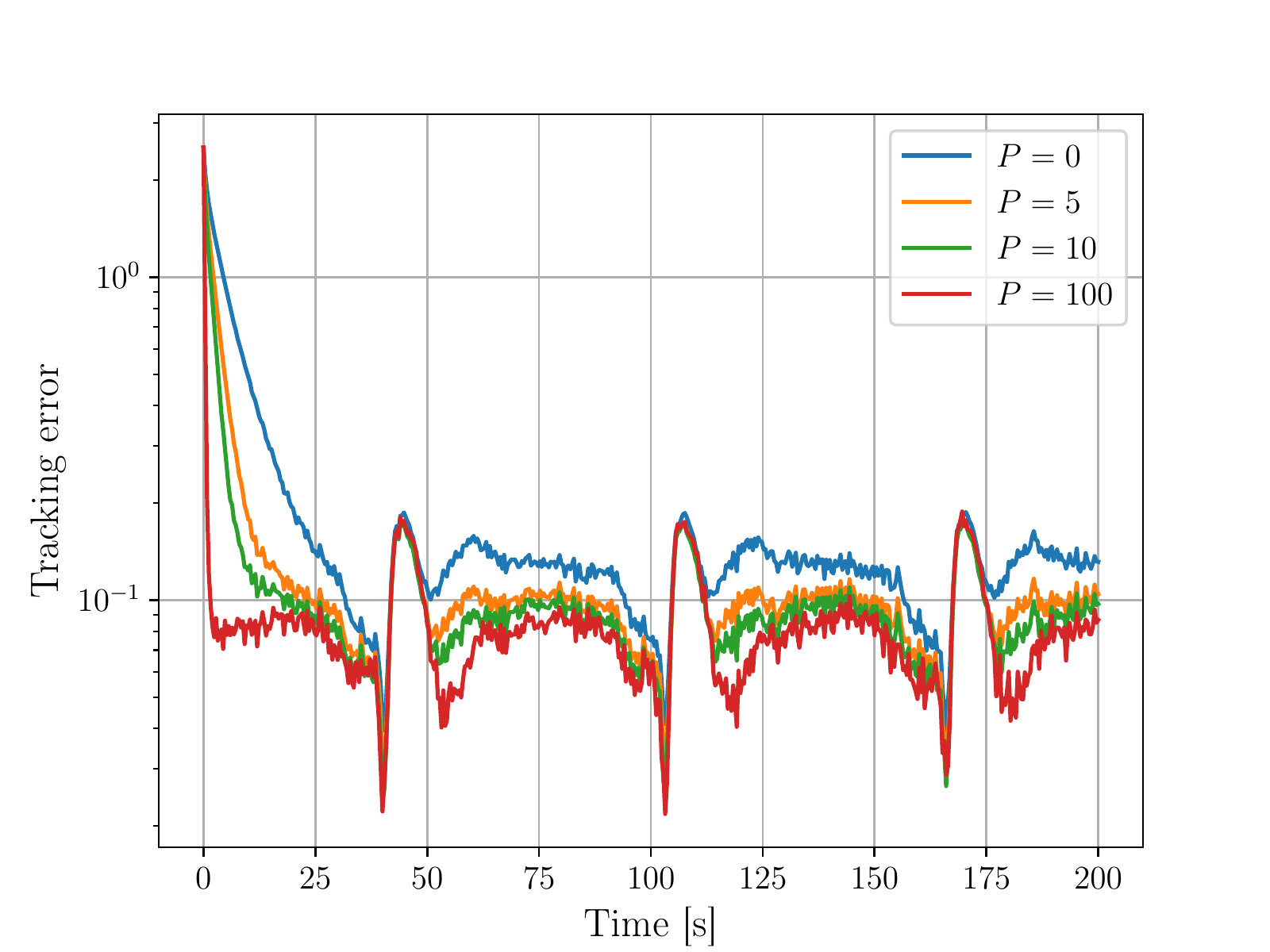}}
  \centerline{(a) Tracking error for different prediction horizons.}\medskip
\end{minipage}
\begin{minipage}[b]{.5\linewidth}
  \centering
  \centerline{\includegraphics[width=1.1\textwidth]{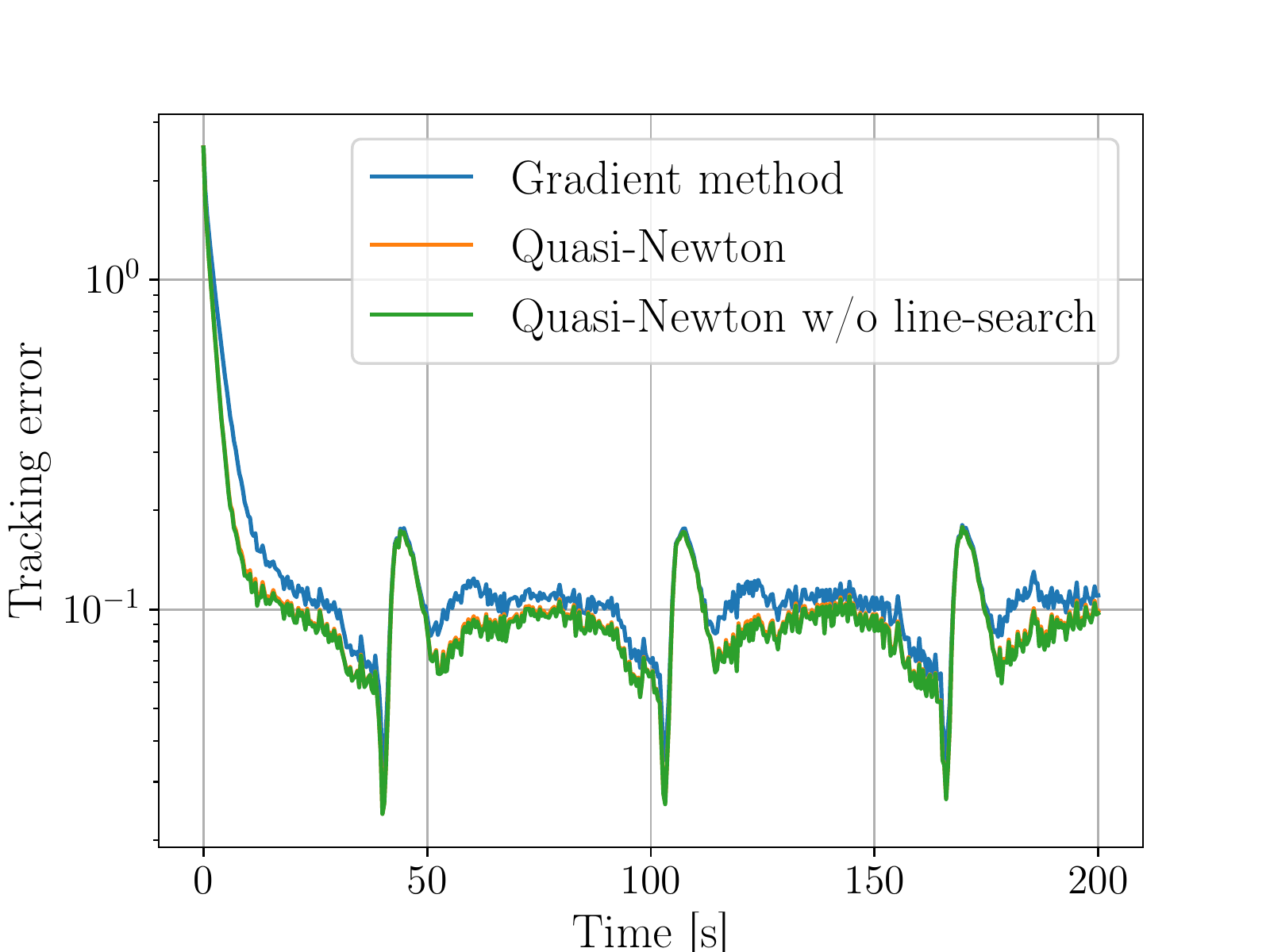}}
  \centerline{(b) Gradient vs quasi-Newton.}\medskip
\end{minipage}
\hfill
\begin{minipage}[b]{0.5\linewidth}
  \centering
  \centerline{\includegraphics[width=1.1\textwidth]{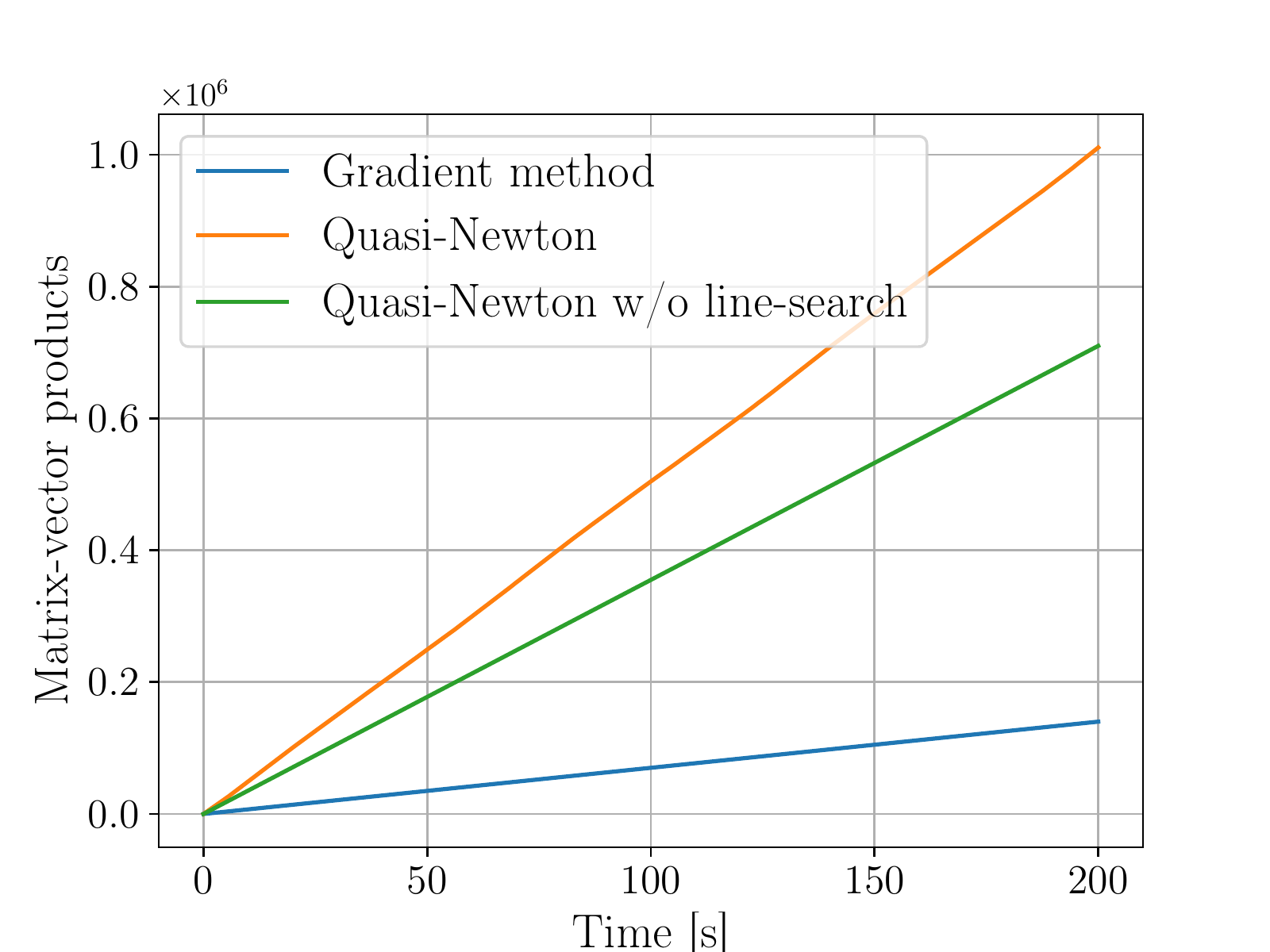}}
  \centerline{(c) Matrix-vector products.}\medskip
\end{minipage}
\caption{Experimental results.}
\label{fig:results}
\end{figure}

In Fig.~\ref{fig:results}(a) we present the evolution of the error 
$$
E_{\textrm{r}} = 	\norm{\x_k - \x_k^*}/6,
$$ 
(the error divided by the number of non-zero components), labeled as ``Tracking error'', for different values of the prediction horizon $P$, and $\alpha = 0.8$, obtained with $T_\mathrm{s} = 0.1 \si{s}$, $C=5$, $\gamma = 0.8/L$ and minimizing the FBE with the line-search quasi-Newton
. Notice that a larger number of prediction steps yields a faster convergence rate and a lower error, which justifies the use of the prediction-correction scheme. Indeed in case we perform only a correction ($P=0$) we obtain the worst performance, which means that the ability to predict the future solution enhances the performance of the optimization algorithm. Note that, even with $P=5$, the performance is better. The quasi-cyclic nature of the error is due to the sinusoidal reference signal.


As noted in Remark~\ref{rem:convexity-fbe}, a convex quadratic cost function $f$ guarantees that the FBE is strongly convex. In this scenario, we can compute a descent direction for the FBE without the need for the line-search procedure included in the quasi-Newton algorithm of \cite{stella2017forward}; moreover, the quasi-Newton algorithms with and without line-search yield the same results. We can also think of applying the gradient method, which has good convergence properties for strongly and well-conditioned convex functions and which does not require the computation of the Hessian of the FBE.

Fig.~\ref{fig:results}(b) depicts the evolution of the error for the quasi-Newton and the gradient methods (with the parameter $P=10, C=5$), while Fig.~\ref{fig:results}(c) the number of matrix-vector products required by the quasi-Newton with and without line-search, and the gradient. Therefore we can choose between a more computationally demanding quasi-Newton, or a gradient method that is simpler to implement but obtains results close to the former.

\section{Conclusion}\label{sec:conclusions}
In this paper, we presented a prediction-correction scheme for time-varying optimization that employs the forward-backward envelope. We described two Theorems that guarantee the convergence of the solution computed by the algorithm to a neighborhood of the optimal solution that depends on the sampling time. Finally we validated the proposed algorithm with some numerical results. Future works will address the problem of relaxing the requirements on the cost functions, especially strong convexity, and perhaps convexity altogether, together with extensive numerical validation and comparison with state-of-the-art methods in specific applications, \emph{e.g.}, dynamic $\ell_1$ reconstruction. 

\addtolength{\textheight}{-2cm}

\addtolength{\textheight}{+2cm}

\appendix

\newpage
\section{Proofs}\label{sec:error-bound}
The sources of error for the proposed algorithm are the Taylor approximation in the prediction step, and the early termination of the quasi-Newton minimization of the FBE.

\subsection{Approximation error}
Let $\afbe(\x)$ be the FBE computed for $h_k$ and $g$, explicitly defined as
\begin{equation*}
	\afbe(\x) = \min_\y \Bigg\{ h_k(\x) +\langle \nabla h_k(\x), \y - \x \rangle  + g(\y) + \frac{\norm{\y - \x}^2}{2\gamma} \Bigg\},
\end{equation*}
and thus that solves the prediction problem~\eqref{eq:approx-gen-equation}\footnote{Notice that since $h_k$ depends only on $\x$ we omit to specify that the gradient is computed w.r.t. $\x$.}. 

\noindent Minimizing this FBE requires that we solve the generalized equation $\nabla \afbe(\x) \ni 0$, and we consider its parametrized version, defined as
\begin{equation}\label{eq:parametrised-eq}
	\nabla \afbe(\x) \ni \p,
\end{equation}
where $\p \in \R^n$.

\noindent Therefore the problem becomes that of finding the solution mapping
$$
	S(\p) = \left\{ \x \in \mathbb{R}^n\ |\ F(\p,\x) = 0 \right\},\ \p\in\mathbb{R}^d
$$
where $F:\mathbb{R}^d \times \mathbb{R}^n \to \mathbb{R}^n$ is defined as $F(\p,\x) = \nabla \afbe(\x) - \p$. Reformulating the problem in this fashion we can make use of Dini's theorem \cite[Th.~1B.1]{dontchev2014implicit}, reported here for convenience.

\begin{theorem}[Dini {\cite[Th.~1B.1]{dontchev2014implicit}}]\label{th:dini}
Let $F:\mathbb{R}^n \times \mathbb{R}^n \to \mathbb{R}^n$ be continuously differentiable in a neighborhood of $(\bar{\p},\bar{\x})$ and such that $F(\bar{\p},\bar{\x}) = 0$ and $\nabla_{\x} F(\bar{\p},\bar{\x})$ is nonsingular. Then the solution mapping $S(\p)$ has a single-valued localization $s$ around $\bar{\p}$ for $\bar{\x}$ which is continuously differentiable in a neighborhood of $\bar{\p}$ with Jacobian satisfying
$$
	\nabla s(\p) = -\nabla_{\x} F(\p,s(\p))^{-1} \nabla_p F(\p,s(\p)).
$$ \endstatement
\end{theorem}

Assume now that $\bar{\x}$ is a solution to~\eqref{eq:parametrised-eq} with $\bar{\p}$ such that $s(\bar{\p}) = \x$. Then Theorem~\ref{th:dini} holds for $F(\p,\x) = \nabla \afbe(\x) - \p$ in a neighborhood of $(\bar{\p},\bar{\x})$ if $F$ is continuously differentiable and has non-singular Jacobian. These conditions are analyzed in the following:

\begin{itemize}
	\item $F$ is differentiable in $\p$ everywhere, while to be differentiable in $\x$ it must be that $\nabla_{\x} \afbe$ is differentiable, which is guaranteed in a neighborhood of $\bar{\x}$ by \cite[Theorem~4.7]{themelis2018forward} under Assumption~\ref{as:first}.
	\item By \cite[Theorem~2.11]{stella2017forward} the Hessian of the FBE in a \textit{strong local minimum}\footnote{A minimum $x^*$ of function $h$ is said to be \textit{locally strong} if there exists $\alpha>0$ such that $h(x)-h(x^*) \geq \alpha\norm{x-x^*}^2$ for any $x$ in a neighborhood of $x^*$.} is positive definite, therefore $\nabla_{\x} F(\bar{\p},\bar{\x}) = \nabla^2 \afbe(\bar{\x})$ is nonsingular if $\bar{\x}$ is a strong local minimum.
\end{itemize}

\noindent Hence as long as $\bar{\x}$ is a strong local minimum it is possible to prove that Theorem~\ref{th:dini} holds for the problem at hand.

The next step is to apply Dini's theorem to provide an upper bound to the error. First notice that a function that is continuously differentiable in a point $\y$ has Lipschitz constant the norm of its gradient in $\y$ \cite[p.~30]{dontchev2014implicit}. Therefore, since $\nabla_p F(\p,s(\p)) = I$, then $s(\p)$ is Lipschitz continuous in a neighborhood of $\bar{\p}$ with constant
$$
	\norm{\nabla s(\bar{\p})} = \norm{\nabla_{\x} F(\bar{\p},\bar{\x})^{-1}} = \norm{\nabla^2 \afbe(\bar{\x})^{-1}} =: K.
$$

As mentioned at the beginning of this proof, the aim is to define an upper bound for the error introduced by the approximation of the cost function, that is $\norm{\bar{\x} - \x^*(t_{k+1})}$. This is accomplished by using the Lipschitz continuity of $s$
$$
	\norm{\bar{\x}-\x^*(t_{k+1})} = \norm{s(\bar{\p})-s(\mathbf{q})} \leq K \norm{\bar{\p}-\mathbf{q}}
$$
where it is necessary to find $\mathbf{q}$ such that $s(\mathbf{q}) = \x^*(t_{k+1})$. Notice that $\x_{k+1}^* := \x^*(t_{k+1})$ must satisfy $\nabla \env(\x_{k+1}^*) \ni 0$ where $\env(\x_{k+1}^*)$ is the forward-backward envelope defined in~\eqref{eq:fbe}. Therefore by choosing
$$
	\mathbf{q} = \nabla \afbe(\x_{k+1}^*) - \nabla \env(\x_{k+1}^*)
$$
we have $s(\mathbf{q}) = \x_{k+1}^*$, since
$$
	\nabla \afbe(\x_{k+1}^*) \ni \mathbf{q} = \nabla \afbe(\x_{k+1}^*) - \nabla \env(x_{k+1}^*)
$$
is verified because by definition $\nabla \env(\x_{k+1}^*) \ni 0$.

Finally it follows that
$$
	\norm{\bar{\x}-\x_{k+1}^*} \leq K \norm{\nabla \afbe(\x_{k+1}^*) - \nabla \env(\x_{k+1}^*)}.
$$

By \cite[Theorem~2.6]{stella2017forward} it follows that
\begin{align}\label{eq:envelope-S-and-R}
	\nabla \env(\x_{k+1}^*) &= (I-\gamma \nabla_{\x\x} f(\x_{k+1}^*;t_{k+1})) R_\gamma(\x_{k+1}^*) \nonumber \\
	&=: S_\gamma(\x_{k+1}^*) R_\gamma(\x_{k+1}^*)
\end{align}
where $R_\gamma(\y) = \gamma^{-1} (\y - \prox_{\gamma g}(\y - \gamma \nabla_{\x} f(\y;t_{k+1})))$ is the \textit{residual}. However, by the definition of residual it holds $R_\gamma(\x_{k+1}^*) = 0$ which actually implies that $\nabla \env(\x_{k+1}^*) = 0$.

\noindent Similarly for $\afbe(\x)$ it holds
$$
	\nabla \afbe(\x_{k+1}^*) = \tilde{S}_\gamma(\x_{k+1}^*) \tilde{R}_\gamma(\x_{k+1}^*)
$$
with $\tilde{S}_\gamma(\x_{k+1}^*)$ and $\tilde{R}_\gamma(\x_{k+1}^*)$ defined substituting $h_k(\cdot)$ to $f(\cdot;t_{k+1})$ in $S_\gamma(\x_{k+1}^*)$ and $R_\gamma(\x_{k+1}^*)$.

Therefore the approximation error can be upper bounded as
\begin{equation}\label{eq:x_bar-x_star}
	\norm{\bar{\x}-\x_{k+1}^*} \leq K \norm{\nabla \afbe(\x_{k+1}^*)}.
\end{equation}

We now proceed to bound $K =\|\nabla^2 \afbe(\bar{\x})^{-1}\|$ and $J := \|\nabla \afbe(\x_{k+1}^*)\|$.

First of all, using the fact that $\nabla^2 h_k(\y) = \nabla_{\x\x} f(\x_k;t_k) =: \Q_k$ it follows $\tilde{S}_\gamma(\x_{k+1}^*) = I - \gamma \Q_k$ which is a positive definite matrix. Indeed observe that by Assumption~\ref{as:first} the maximum and minimum eigenvalues of $\Q_k$ are, respectively: $\lambda_M(\Q_k) = L$ and $\lambda_m(\Q_k) = m$, with $L \geq m > 0$. Therefore since the eigenvalues of $\tilde{S}_\gamma(\x)$ are of the form $\lambda_i(\tilde{S}_\gamma(\x)) = 1 - \gamma \lambda_i(\Q_k)$ it follows that
\begin{equation}\label{eq:eigenvalues-S}
	\lambda_M(\tilde{S}_\gamma(\x)) \leq 1-\gamma m \quad\text{and}\quad \lambda_m(\tilde{S}_\gamma(\x)) \geq 1-\gamma L,
\end{equation}
with $1-\gamma L > 0$ because $\gamma < 1/L$.

\noindent Furthermore it holds by \cite[Theorem~2.10]{stella2017forward}
$$
	\nabla^2 \afbe(\bar{\x}) = \gamma^{-1} \tilde{S}_\gamma(\bar{\x}) (I-\tilde{P}_\gamma(\bar{\x})\tilde{S}_\gamma(\bar{\x}))
$$
where $\tilde{P}_\gamma(\bar{\x})$ is semi-definite positive and $\|\tilde{P}_\gamma(\x)\| \leq 1$.

\subsubsection{Computation of $J$}
By the formula for the gradient of $\afbe(\x)$ it holds
\begin{equation}\label{eq:norm-afbe}
	\|\nabla \afbe(\x_{k+1}^*)\| \leq \|\tilde{S}_\gamma(\x_{k+1}^*)\|\|\tilde{R}_\gamma(\x_{k+1}^*)\|
\end{equation}
where by~\eqref{eq:eigenvalues-S} it holds $\|\tilde{S}_\gamma(\x_{k+1}^*)\| \leq 1 - \gamma m$ . Therefore to compute $J$ it is necessary to provide an upper bound to the norm of the residual $\|\tilde{R}_\gamma(\x_{k+1}^*)\|$.

Recalling the definition of the residual and by the fact that $R_\gamma(\x_{k+1}^*) = 0$ it holds
\begin{align*}
	& \|\tilde{R}_\gamma(\x_{k+1}^*)\| = \|R_\gamma(\x_{k+1}^*) - \tilde{R}_\gamma(\x_{k+1}^*)\| \\
	& \qquad =\gamma^{-1} \| \prox_{\gamma g}(\x_{k+1}^* - \gamma \nabla_{\x} f(\x_{k+1}^*;t_{k+1})) \\
	&\qquad\qquad\qquad\qquad\qquad - \prox_{\gamma g} (\x_{k+1}^* - \gamma \nabla h_k(\x_{k+1}^*)) \| \\
	& \qquad \leq \norm{\nabla_{\x} f(\x_{k+1}^*;t_{k+1}) - \nabla h_k(\x_{k+1}^*))}
\end{align*}
where the nonexpansiveness of the proximal operator was used to derive the last inequality.

By using the definition of $h_k$, it follows
\begin{align*}
	& \norm{\nabla_{\x} f(\x_{k+1}^*;t_{k+1}) - \nabla h_k(\x_{k+1}^*))} \\
	& \leq \norm{\nabla_{\x} f(\x_{k+1}^*;t_{k+1}) - \nabla_{\x} f(\x_k;t_k)} + \\
	& \qquad\qquad\qquad + \norm{\Q_k} \norm{\x_{k+1}^* - \x_k} + T_\mathrm{s} \norm{\nabla_{t\x} f(\x_k;t_k)} \\
	& \leq \norm{\nabla_{\x} f(\x_{k+1}^*;t_{k+1}) - \nabla_{\x} f(\x_k;t_k)} + L \norm{\x_{k+1}^* - \x_k} + T_\mathrm{s} C_0
\end{align*}
where the upper bounds in Assumptions~\ref{as:first}-\ref{as:first-bis} on the derivatives of $f$ were used.

The first term on the right-hand side of the inequality remains now to be computed. It holds
\begin{align*}
	& \norm{\nabla_{\x} f(\x_{k+1}^*;t_{k+1}) - \nabla_{\x} f(\x_k;t_k)} \\
	& \leq \norm{\nabla_{\x} f(\x_{k+1}^*;t_{k+1}) - \nabla_{\x} f(\x_k;t_{k+1})} + \\
	& \qquad\qquad\qquad\qquad + \norm{\nabla_{\x} f(\x_k;t_{k+1}) - \nabla_{\x} f(\x_k;t_k)}
\end{align*}
where by Lipschitz continuity
$$
	\norm{\nabla_{\x} f(\x_{k+1}^*;t_{k+1}) - \nabla_{\x} f(\x_k;t_{k+1})} \leq L \norm{\x_{k+1}^* - \x_k}
$$
and $\norm{\nabla_{\x} f(\x_k;t_{k+1}) - \nabla_{\x} f(\x_k;t_k)} \leq T_\mathrm{s}C_0$ (see \cite[Appendix~A]{simonetto2017prediction}).

Finally the results above yield
\begin{equation}\label{eq:bound-residual}
	\|\tilde{R}_\gamma(\x_{k+1}^*)\| \leq 2(L \norm{\x_{k+1}^* - \x_k} + T_\mathrm{s}C_0)
\end{equation}
and therefore
$$
	\norm{\bar{\x} - \x_{k+1}^*} \leq 2 K (1-\gamma m) (L \norm{\x_k - \x_{k+1}^*} + T_\mathrm{s}C_0).
$$

The problem now is to bound $\norm{\x_k - \x_{k+1}^*}$ in terms of $\norm{\x_k - \x_k^*}$; first of all, it holds $\norm{\x_k - \x_{k+1}^*} \leq \norm{\x_{k+1}^* - \x_k^*} + \norm{\x_k - \x_k^*}$ and so an upper bound for $\norm{\x_{k+1}^* - \x_k^*}$ must be found.

Recall that with $\x^*(t)$ we denote the optimal solution of the original time-varying problem~\eqref{eq:tv-problem}. Notice that Dini's theorem holds for the corresponding FBE at time $t$: $\nabla_{\x} \env(\x;t)$ with $t \in \mathbb{R}_+$ around $\x^*(t)$. Therefore the solution mapping $\x^*(t)$ has a Lipschitz constant upper bounded by
\begin{align*}
	& \norm{\nabla_{\x\x} \env(\x^*(t);t)^{-1}} \norm{\nabla_{t\x} \env(\x^*(t);t)} \\ & \qquad\qquad =:  \bar{K} \norm{\nabla_{t\x} \env(\x^*(t);t)}.
\end{align*}
The term $\bar{K}$ will be computed in the next section alongside $K$, while the second term requires evaluating the time derivative of the gradient of the FBE. In particular, by derivative rules, it holds that
$$
	\nabla_{t\x} \env(\x;t) = \nabla_t S_\gamma(\x;t) R_\gamma(\x;t) + S_\gamma(\x;t) \nabla_t R_\gamma(\x;t)
$$
where the first term on the right-hand side can be ignored, since at $(\x^*(t);t)$ it is zero. By using the definition~\eqref{eq:envelope-S-and-R}, it holds that $S_\gamma(\x;t) = I-\gamma\nabla_{\x\x}f(\x;t)$, whose norm is upper bounded by $1-\gamma m$; hence it remains only to compute an upper bound to the norm of $\nabla_t R_\gamma(\x;t)$.

By computing the time derivative of the residual, it follows that
\begin{align*}
	\nabla_t R_\gamma(\x;t) &= - \gamma^{-1} \nabla_t [\prox_{\gamma g}(\x - \gamma \nabla_{\x} f(\x;t))] \\
	&= \nabla_{t\x} f(\x;t) J\prox_{\gamma g}(\x - \gamma\nabla_{\x} f(\x;t)) \\
	&= \nabla_{t\x} f(\x;t) P_\gamma(\x;t)
\end{align*}
where the second inequality is derived using the chain rule, and the third by using the definition of the semi-definite positive matrix $P_\gamma(\x;t)$ reported in \cite{stella2017forward}.

Therefore it holds that
$$
	\norm{\nabla_t R_\gamma(\x^*(t);t)} \leq \norm{\nabla_{t\x} f(\x^*(t);t)} \norm{P_\gamma(\x^*(t);t)} \leq C_0
$$
since $\norm{P_\gamma(\x^*(t);t)} \leq 1$.

Finally, the solution mapping $\x^*(t)$ has Lipschitz constant $\bar{K}(1-\gamma m)C_0$ and so
\begin{align}\label{eq:difference-optima}
	& \norm{\x_{k+1}^* - \x^*_k} = \norm{\x^*_{k+1} - \x^*_k} \\
	& \qquad \leq \bar{K}(1-\gamma m)C_0 |t_{k+1} - t_k| = \bar{K}(1-\gamma m)C_0T_\mathrm{s} \nonumber.
\end{align}

\subsubsection{Computation of $K$}
Let $A,B \in \mathbb{R}^{n \times n}$ with $A$ nonsingular, and denote with $\svmin{\cdot}$ and $\svmax{\cdot}$ the minimum and maximum singular values of a matrix, respectively. The following facts hold true for symmetric nonsingular matrices\cite[Ch.~5]{meyer2000matrix}.

\begin{enumerate}
	\item $\norm{A^{-1}} = \svmax{A^{-1}} = 1/\svmin{A}$;
	\item $\svmin{AB} \geq \svmin{A} \svmin{B}$.
\end{enumerate}

\noindent Therefore it follows
$$
	K = \|\nabla^2 \afbe(\bar{\x})^{-1}\| = 1/\svmin{\nabla^2 \afbe(\bar{\x})}
$$
which can be bounded using lower bounds for $\svmin{\nabla^2 \afbe(\bar{\x})}$. In particular,
\begin{equation}\label{eq:svmin_afbe}
	\svmin{\nabla^2 \afbe(\bar{\x})} \geq \gamma^{-1} \svmin{\tilde{S}_\gamma(\bar{\x})} \svmin{(I-\tilde{P}_\gamma(\bar{\x})\tilde{S}_\gamma(\bar{\x}))}
\end{equation}
and the problem is to compute the minimum singular values of $\tilde{S}_\gamma(\bar{\x})$ and $(I-\tilde{P}_\gamma(\bar{\x})\tilde{S}_\gamma(\bar{\x}))$.

First of all, for a symmetric matrix $A \in \mathbb{R}^{n \times n}$ it holds $\svmin{A} = |\lambda_m(A)|$. Recalling~\eqref{eq:eigenvalues-S} we have that $\svmin{\tilde{S}_\gamma(\bar{\x})} \geq 1-\gamma L > 0$. Since $I - \tilde{P}_\gamma(\bar{\x})\tilde{S}_\gamma(\bar{\x})$ is symmetric and positive definite, then
$$
	\svmin{(I-\tilde{P}_\gamma(\bar{\x})\tilde{S}_\gamma(\bar{\x}))} = \lambda_m\left((I-\tilde{P}_\gamma(\bar{\x})\tilde{S}_\gamma(\bar{\x}))\right)
$$
and it is necessary to lower bound the spectrum of $(I-\tilde{P}_\gamma(\bar{\x})\tilde{S}_\gamma(\bar{\x}))$.

It holds that 
$$
\lambda_m\left((I-\tilde{P}_\gamma(\bar{\x})\tilde{S}_\gamma(\bar{\x}))\right) = 1 - \lambda_M\left(\tilde{P}_\gamma(\bar{\x})\tilde{S}_\gamma(\bar{\x})\right)
$$ 
and an upper bound for $\lambda_M\left(\tilde{P}_\gamma(\bar{\x})\tilde{S}_\gamma(\bar{\x})\right)$ must be found. Notice that $\tilde{P}_\gamma(\bar{\x})\tilde{S}_\gamma(\bar{\x})$ is symmetric and positive definite \cite[Appendix~B]{stella2017forward}. The following inequalities hold
\begin{align*}
	& \lambda_M\left(\tilde{P}_\gamma(\bar{\x})\tilde{S}_\gamma(\bar{\x})\right) = \norm{\tilde{P}_\gamma(\bar{\x})\tilde{S}_\gamma(\bar{\x})} \\
	& \qquad\qquad \leq \norm{\tilde{P}_\gamma(\bar{\x})} \norm{\tilde{S}_\gamma(\bar{\x})} \leq \norm{\tilde{S}_\gamma(\bar{\x})} \leq 1-\gamma m
\end{align*}
where the fact that $\norm{\tilde{P}_\gamma(\bar{\x})} \leq 1$ and the result~\eqref{eq:eigenvalues-S} were used.

Finally, these results yield
\begin{align}\label{eq:svmin_I_PS}
	& \svmin{(I-\tilde{P}_\gamma(\bar{\x})\tilde{S}_\gamma(\bar{\x}))} = \lambda_m\left((I-\tilde{P}_\gamma(\bar{\x})\tilde{S}_\gamma(\bar{\x}))\right) \nonumber \\
	& \geq 1 - \lambda_M\left(\tilde{P}_\gamma(\bar{\x})\tilde{S}_\gamma(\bar{\x})\right) \geq 1 - (1 - \gamma m) = \gamma m.
\end{align}

\noindent Finally, substituting~\eqref{eq:svmin_I_PS} and $\svmin{\tilde{S}_\gamma({\bar{\x})}} \geq 1-\gamma L$ into~\eqref{eq:svmin_afbe}, it follows
$$
	\svmin{\nabla^2 \afbe(\bar{\x})} \geq \gamma^{-1} (1-\gamma L) \gamma m = m (1- \gamma L)
$$
and thus
\begin{equation}\label{eq:bound-K}
	K \leq \frac{1}{m(1-\gamma L)}.
\end{equation}
Notice that the same bound holds for $\bar{K}$ as well since, like $f$, it is $h_k \in \mathcal{S}_{m,L}(\R^n)$.

\subsubsection{Approximation error bound}
The results derived in the previous sections can be used to bound the approximation error as follows
\begin{align}
	\norm{\bar{\x} - \x_{k+1}^*} &\leq \frac{2(1-\gamma m)}{m (1-\gamma L)} \Bigg(L\norm{\x_k - \x^*_k} \nonumber \\
	& \qquad\qquad + \frac{L(1-\gamma m)}{m(1-\gamma L)}T_\mathrm{s}C_0 + T_\mathrm{s} C_0\Bigg) \nonumber \\
	&\leq a_1 \norm{\x_k - \x^*_k} + a_0 \label{eq:approximation-error}
\end{align}
where the coefficients are defined as
\begin{align*}
	a_1 &= \frac{2L(1-\gamma m)}{m(1-\gamma L)} \\
	a_0 &= 2C_0T_\mathrm{s} \frac{1-\gamma m}{m(1-\gamma L)} \left[\frac{L(1-\gamma m)}{m(1-\gamma L)} + 1\right]
\end{align*}
with $a_0$ that linearly depends on the sampling time $T_\mathrm{s}$.

\subsection{Early termination errors}
The previous section derived an upper bound to the approximation error introduced by the approximate FBE used in the prediction step instead of the correct FBE. However there are two other sources of errors, namely the early termination of the minimisation algorithms applied to the $\afbe(\x)$ and the FBE in the prediction and correction steps, respectively. These will be the focus of the current section.

Assume that the minimisations are carried out using the \textit{quasi-Newton} Algorithm~2 in \cite{stella2017forward} applied to the $\afbe$ or the $\env$. This method is described in Algorithm \ref{alg:quasi-newton-fbe} in the general case.

\begin{algorithm}
\caption{Quasi-Newton method for the FBE.}
\label{alg:quasi-newton-fbe}
\begin{algorithmic}[1]
	\Require $\hat{\x}^0 \in \mathbb{R}^n$, envelope parameter $\gamma \in (0, 1/L)$, maximum number of iterations $I$.
	\State $i \gets 0$
	\While{$i < I$}
		\State compute the nonsingular matrix $B^i$ using the BFGS method
		\State compute the descent direction
		$$
			\mathbf{d}^i = -(B^i)^{-1}\nabla \env(\hat{\x}^i)
		$$
		\State select the step-size $\tau^i \geq 0$ s.t. $\env(\w^i) \leq \env(\hat{\x}^i)$ where $\w^i = \hat{\x}^i + \tau^i \mathbf{d}^i$
		\State $\hat{\x}^{i+1} \gets \prox_{\gamma g}(\x - \gamma\nabla_{\x} f(\x;t_k))$
		\State $i \gets i+1$
	\EndWhile
\end{algorithmic}
\end{algorithm}

\begin{remark}
In case $L$ is not known, it is possible to apply Algorithm~1 of \cite{stella2017forward} which includes a line-search procedure for choosing $\gamma$.
\end{remark}

The sequence $\{\hat{\x}^m\}_{m\in\mathbb{N}}$ produced by Algorithm \ref{alg:quasi-newton-fbe} can be proved to converge to a critical point with a \textit{super-linear} rate, by a combination of Theorems 3.6 and 4.3 in \cite{stella2017forward}. Moreover the result holds globally by strong convexity.

Since super-linear convergence implies linear convergence as well, then there exist $\zeta \in (0,1)$ (see~\eqref{eq:zeta} for an estimate of its value) such that
$$
	\norm{\hat{\x}^{p+1} - \bar{\x}} \leq \zeta \norm{\hat{\x}^p - \bar{\x}} \quad p=1,2,\ldots,P
$$
and
$$
	\norm{\hat{\x}^{c+1} - \x_{k+1}^*} \leq \zeta \norm{\hat{\x}^c - \x_{k+1}^*} \quad c=1,2,\ldots,C
$$
where $\hat{\x}^p$ and $\hat{\x}^c$ are the dummy variables used during prediction and correction.

The dummy variables are initialized as follows: for the prediction $\hat{\x}^0 = \x_k$, for the correction $\hat{\x}^0 = \tilde{\x}_{k+1|k} = \hat{\x}^P$. Therefore iterating
\begin{align*}
	\norm{\tilde{\x}_{k+1|k} - \bar{\x}} = \norm{\tilde{\x}_{k+1|k} - \x_{k+1|k}} &\leq \zeta^P \norm{\x_k - \x_{k+1|k}} \\
	\norm{\x_{k+1} - \x_{k+1}^*} &\leq \zeta^C \norm{\tilde{\x}_{k+1|k} - \x_{k+1}^*}.
\end{align*}

\subsection{Overall error bound}
During the previous sections the following bounds have been derived
\begin{align}
	\norm{\bar{\x} - \x_{k+1}^*} &\leq a_1 \norm{\x_k - \x^*(t_k)} + a_0 \label{eq:ineq-1} \\
	\norm{\tilde{\x}_{k+1|k} - \bar{\x}} &\leq \zeta^P \norm{\x_k - \x_{k+1|k}} \label{eq:ineq-2} \\
	\norm{\x_{k+1} - \x_{k+1}^*} &\leq \zeta^C \norm{\tilde{\x}_{k+1|k} - \x_{k+1}^*} \label{eq:ineq-3}
\end{align}
and therefore the last thing to do is to combine them to derive a bound for the error $\norm{\x_{k+1} - \x_{k+1}^*}$. Following the same steps of \cite[Appendix~B]{simonetto2017prediction} from inequalities \eqref{eq:ineq-1}, \eqref{eq:ineq-2} and the results above it is possible to compute
\begin{align*}
	&\norm{\tilde{\x}_{k+1|k} - \x_{k+1}^*} \leq \norm{\tilde{\x}_{k+1|k} - \bar{\x}} + \norm{\bar{\x} - \x_{k+1}^*} \\
	&\leq \zeta^P \norm{\x_k - \bar{\x}} + \norm{\bar{\x} - \x_{k+1}^*} \\
	&\leq \zeta^P (\norm{\x_k - \x_k^*} + \norm{\x_k^* - \x_{k+1}^*} + \norm{\bar{\x} - \x_{k+1}^*}) + \\&\qquad\qquad\qquad\qquad\qquad\qquad\norm{\bar{\x} - \x_{k+1}^*} \\
	&\leq \zeta^P \norm{\x_k - \x_k^*} + \zeta^P \frac{(1-\gamma m)T_\mathrm{s}C_0}{m(1-\gamma L)} +\\ & \qquad\qquad\qquad\qquad\qquad\qquad (\zeta^P + 1) \norm{\bar{\x} - \x_{k+1}^*} \\
	&\leq [\zeta^P + a_1 (\zeta^P + 1)] \norm{\x_k - \x_k^*} \\ & \qquad\qquad + \zeta^P \frac{(1-\gamma m)T_\mathrm{s}C_0}{m(1-\gamma L)} + (\zeta^P + 1) a_0.
\end{align*}
Using now inequality \eqref{eq:ineq-3} it follows
\begin{equation}
	\norm{\x_{k+1} - \x_{k+1}^*} \leq A_1 \norm{\x_k - \x_k^*} + A_0
\end{equation}
where
\begin{align*}
	A_1 &= \zeta^C [\zeta^P + a_1(\zeta^P + 1)] \\
	&= \zeta^C \left[\zeta^P + (\zeta^P+1) \frac{2L(1-\gamma m)}{m(1-\gamma L)}\right] \\
	A_0 &= \zeta^C \left[\zeta^P \frac{(1-\gamma m)T_\mathrm{s}C_0}{m(1-\gamma L)} + (\zeta^P + 1)a_0 \right].
\end{align*}

Therefore, for the algorithm to converge to a bounded error it is necessary that $A_1 < 1$, which must be guaranteed by choosing suitable prediction and correction horizons. Notice that $A_0$ depends linearly on the sampling time $T_\mathrm{s}$, therefore it holds that
$$
	\limsup_{k\to\infty} \norm{\x_{k+1} - \x_{k+1}^*} = O(\zeta^C T_\mathrm{s}).
$$
We have thus proved Theorem~\ref{th:linear-convergence}. \hfill$\blacksquare$

\subsubsection{$O(T^2_\mathrm{s})$ convergence}
Suppose now that Assumption~\ref{as:second} holds, it is possible to show that the error bound tends to $O(T^2_\mathrm{s})$.

Consider that
\begin{equation}\label{eq:R-epsilon}
	\norm{\tilde{R}_\gamma(\x_{k+1}^*)} \leq \norm{\nabla_{\x} f(\x_{k+1}^*;t_{k+1}) - \nabla h_k(\x_{k+1}^*)} \leq \norm{\pmb{\epsilon}},
\end{equation}
where $\pmb{\epsilon}$ is the residual of the Taylor expansion of $\nabla_{\x} f(\cdot;t_{k+1})$. Under Assumption~\ref{as:second}, the Taylor residual can be bound as follows
\begin{align}
	\norm{\pmb{\epsilon}} &\leq \frac{1}{2}\Big( \norm{\nabla_{\x\x\x}f(\x_k;t_k)} \norm{\x_{k+1}^*-\x_k}^2 + \nonumber \\
	& + T_\mathrm{s} \norm{\nabla_{t\x\x}f(\x_k;t_k)} \norm{\x_{k+1}^*-\x_k} + \nonumber \\
	& + T_\mathrm{s} \norm{\nabla_{\x t\x}f(\x_k;t_k)} \norm{\x_{k+1}^*-\x_k} +\nonumber\\&\qquad\qquad\qquad\qquad\qquad\qquad +T_\mathrm{s}^2 \norm{\nabla_{tt\x}f(\x_k;t_k)} \Big) \nonumber \\
	&\leq \frac{C_1}{2}\norm{\x_{k+1}^*-\x_k}^2 + T_\mathrm{s}C_2\norm{\x_{k+1}^*-\x_k} + \frac{1}{2}T_\mathrm{s}^2C_3. \label{eq:epsilon-bound}
\end{align}

We want now to use this result to compute an upper bound to $\norm{\bar{\x} - \x_{k+1}^*}$ that is stricter than~\eqref{eq:ineq-1}. First of all, substituting~\eqref{eq:R-epsilon} into the bound~\eqref{eq:norm-afbe} we obtain
\begin{equation}\label{eq:afbe-epsilon}
	\|\nabla \afbe(\x_{k+1}^*)\| \leq \norm{\tilde{S}_\gamma(\x_{k+1}^*)} \norm{\tilde{R}_\gamma(\x_{k+1}^*)} \leq (1-\gamma m) \norm{\pmb{\epsilon}}
\end{equation}
where we used the bound $\|\tilde{S}_\gamma(\x_{k+1}^*)\| \leq 1-\gamma m$. Therefore applying the bound~\eqref{eq:afbe-epsilon} into~\eqref{eq:x_bar-x_star} yields
\begin{align*}
	\norm{\bar{\x}-\x_{k+1}^*} &\leq K \norm{\nabla \afbe(\x_{k+1}^*)} \\
	&\leq K (1-\gamma m) \norm{\pmb{\epsilon}} \leq \frac{1-\gamma m}{m(1-\gamma L)} \norm{\pmb{\epsilon}}
\end{align*}
where the bound~\eqref{eq:bound-K} for $K$ was used. Substituting~\eqref{eq:epsilon-bound} finally we get 
\begin{equation}\label{eq:x_bar-with-a-coefficients}
	\norm{\bar{\x}-\x_{k+1}^*} \leq a_2 \norm{\x_k-\x_k^*}^2 + a_1 \norm{\x_k-\x_k^*} + a_0
\end{equation}
where
\begin{align*}
	a_2 &= \frac{1-\gamma m}{m(1-\gamma L)}\frac{C_1}{2} \\
	a_1 &= T_\mathrm{s} \frac{1-\gamma m}{m(1-\gamma L)} \left[\frac{1-\gamma m}{m(1-\gamma L)}C_0C_1 + C_2\right] \\
	a_0 &= T_\mathrm{s}^2 \frac{1-\gamma m}{m(1-\gamma L)} \Bigg[ \left(\frac{1-\gamma m}{m(1-\gamma L)}\right)^2 \frac{C_1C_0^2}{2} + \\ & + \frac{1-\gamma m}{m(1-\gamma L)} C_0C_2 + \frac{1}{2}C_3 \Bigg].
\end{align*}
Notice that $a_1$ and $a_0$ linearly depend on $T_\mathrm{s}$ and $T_\mathrm{s}^2$, respectively.

With computations very similar to those carried out during the previous section, exchanging~\eqref{eq:x_bar-with-a-coefficients} for~\eqref{eq:ineq-1}, it is then possible to compute the bound
\begin{equation}
	\norm{\x_{k+1}-\x_{k+1}^*} \leq A_2\norm{\x_k-\x_k^*}^2 + A_1\norm{\x_k-\x_k^*} + A_0
\end{equation}
where
\begin{align*}
	A_2 &= \zeta^C (\zeta^P+1) a_2 \\
	A_1 &= \zeta^C [\zeta^P + a_1 (\zeta^P+1)] \\
	A_0 &= \zeta^C \left[\zeta^P \frac{1-\gamma m}{m(1-\gamma L)}C_0T_\mathrm{s} + (\zeta^P+1)a_0 \right].
\end{align*}

In order to prove convergence, it is now possible to use the argument presented in \cite[Appendix~B]{simonetto2017prediction}, which guarantees convergence if $\tau > \zeta^P \zeta^C$,
\begin{align*}
	& T_\mathrm{s} < \frac{\tau - \zeta^P \zeta^C}{\zeta^C(\zeta^P + 1)} \times \\
	&\qquad \times \left\{ \frac{1-\gamma m}{m(1-\gamma L)} \left[\frac{1-\gamma m}{m(1-\gamma L)}C_0C_1 + C_2\right] \right\}^{-1} =: \bar{T}_\mathrm{s},
\end{align*}
and
$$
	\norm{\x_0 - \x_0^*} \leq \frac{\tau - A_1}{A_2} =: \bar{R}.
$$

\noindent Therefore the asymptotic error satisfies
$$
	\limsup_{k\to\infty} \norm{\x_k - \x_k^*} \leq O(T_\mathrm{s}^2 \zeta^C) + O(T_\mathrm{s} \zeta^C\zeta^P),
$$
which proves Theorem~\ref{th:quadratic-convergence}. \hfill$\blacksquare$


\begin{thebibliography}{10}

\bibitem{asif2010dynamic}
M~Salman Asif and Justin Romberg,
\newblock ``Dynamic updating for $\ell_1$ minimization,''
\newblock {\em IEEE Journal of selected topics in signal processing}, vol. 4,
  no. 2, pp. 421--434, 2010.

\bibitem{asif2014sparse}
M~Salman Asif and Justin Romberg,
\newblock ``Sparse recovery of streaming signals using $\ell_1$-homotopy,''
\newblock {\em IEEE Transactions on Signal Processing}, vol. 62, no. 16, pp.
  4209--4223, 2014.

\bibitem{Vaswani2015}
N.~Vaswani and J.~Zhan,
\newblock ``Recursive recovery of sparse signal sequences from compressive
  measurements: A review,''
\newblock {\em IEEE Transactions on Signal Processing}, vol. 64, no. 13, pp.
  3523 -- 3549, 2016.

\bibitem{Yang2015}
Y.~Yang, M.~Zhang, M.~Pesavento, and D.~P. Palomar,
\newblock ``An online parallel and distributed algorithm for recursive
  estimation of sparse signals,''
\newblock {\em IEEE Transactions on Signal and Information Processing over
  Networks}, vol. 2, no. 3, pp. 290 -- 305, 2016.

\bibitem{charles2016dynamic}
Adam~S Charles, Aurele Balavoine, and Christopher~J Rozell,
\newblock ``Dynamic filtering of time-varying sparse signals via $\ell_1$
  minimization,''
\newblock {\em IEEE Transactions on Signal Processing}, vol. 64, no. 21, pp.
  5644--5656, 2016.

\bibitem{Sopasakis2016}
P.~Sopasakis, N.~Freris, and P.~Patrinos,
\newblock ``Accelerated reconstruction of a compressively sampled data
  stream,''
\newblock in {\em Proceedings of the 24th EUSIPCO}, Budapest, Hungary,
  September 2016, pp. 1078 -- 1082.

\bibitem{jerez2014embedded}
Juan~L Jerez, Paul~J Goulart, Stefan Richter, George~A Constantinides, Eric~C
  Kerrigan, and Manfred Morari,
\newblock ``Embedded online optimization for model predictive control at
  megahertz rates,''
\newblock {\em IEEE Transactions on Automatic Control}, vol. 59, no. 12, pp.
  3238--3251, 2014.

\bibitem{hours2016parametric}
Jean-Hubert Hours and Colin~N Jones,
\newblock ``A parametric nonconvex decomposition algorithm for real-time and
  distributed {NMPC},''
\newblock {\em IEEE Transactions on Automatic Control}, vol. 61, no. 2, pp.
  287--302, 2016.

\bibitem{gutjahr2017lateral}
Benjamin Gutjahr, Lutz Gr{\"o}ll, and Moritz Werling,
\newblock ``Lateral vehicle trajectory optimization using constrained linear
  time-varying {MPC},''
\newblock {\em IEEE Transactions on Intelligent Transportation Systems}, vol.
  18, no. 6, pp. 1586--1595, 2017.

\bibitem{verscheure2009time}
Diederik Verscheure, Bram Demeulenaere, Jan Swevers, Joris De~Schutter, and
  Moritz Diehl,
\newblock ``Time-optimal path tracking for robots: A convex optimization
  approach,''
\newblock {\em IEEE Transactions on Automatic Control}, vol. 54, no. 10, pp.
  2318--2327, 2009.

\bibitem{ardeshiri2011convex}
Tohid Ardeshiri, Mikael Norrl{\"o}f, Johan L{\"o}fberg, and Anders Hansson,
\newblock ``Convex optimization approach for time-optimal path tracking of
  robots with speed dependent constraints,''
\newblock {\em IFAC Proceedings Volumes}, vol. 44, no. 1, pp. 14648--14653,
  2011.

\bibitem{dixit2018online}
Rishabh Dixit, Amrit~Singh Bedi, Ruchi Tripathi, and Ketan Rajawat,
\newblock ``Online learning with inexact proximal online gradient descent
  algorithms,''
\newblock {\em arXiv preprint arXiv:1806.00202}, 2018.

\bibitem{simonetto2016class}
Andrea Simonetto, Aryan Mokhtari, Alec Koppel, Geert Leus, and Alejandro
  Ribeiro,
\newblock ``A class of prediction-correction methods for time-varying convex
  optimization.,''
\newblock {\em IEEE Trans. Signal Processing}, vol. 64, no. 17, pp. 4576--4591,
  2016.

\bibitem{simonetto2017prediction}
Andrea Simonetto and Emiliano Dall’Anese,
\newblock ``Prediction-correction algorithms for time-varying constrained
  optimization,''
\newblock {\em IEEE Transactions on Signal Processing}, vol. 65, no. 20, pp.
  5481--5494, 2017.

\bibitem{simonetto2018tac}
A.~Simonetto,
\newblock ``Dual prediction-correction methods for linearly constrained
  time-varying convex programs,''
\newblock {\em IEEE Transactions on Automatic Control (to appear)}, 2018.

\bibitem{fazlyab2017prediction}
Mahyar Fazlyab, Santiago Paternain, Victor~M. Preciado, and Alejandro Ribeiro,
\newblock ``Prediction-correction interior-point method for time-varying convex
  optimization,''
\newblock {\em IEEE Transactions on Automatic Control}, 2017.

\bibitem{Rahili2015}
S.~Rahili and W.~Ren,
\newblock ``Distributed convex optimization for continuous-time dynamics with
  time-varying cost functions,''
\newblock {\em IEEE Transactions on Automatic Control}, vol. 62, no. 4, pp.
  1590 -- 1605, 2017.

\bibitem{fazlyab2016self}
Mahyar Fazlyab, Cameron Nowzari, George~J. Pappas, Alejandro Ribeiro, and
  Victor~M. Preciado,
\newblock ``Self-triggered time-varying convex optimization,''
\newblock in {\em 2016 IEEE 55th Conference on Decision and Control (CDC)}.
  IEEE, 2016, pp. 3090--3097.

\bibitem{patrinos2013proximal}
Panagiotis Patrinos and Alberto Bemporad,
\newblock ``Proximal newton methods for convex composite optimization,''
\newblock in {\em Decision and Control (CDC), 2013 IEEE 52nd Annual Conference
  on}. IEEE, 2013, pp. 2358--2363.

\bibitem{stella2017forward}
Lorenzo Stella, Andreas Themelis, and Panagiotis Patrinos,
\newblock ``Forward--backward quasi-newton methods for nonsmooth optimization
  problems,''
\newblock {\em Computational Optimization and Applications}, vol. 67, no. 3,
  pp. 443--487, 2017.

\bibitem{themelis2018forward}
Andreas Themelis, Lorenzo Stella, and Panagiotis Patrinos,
\newblock ``Forward-backward envelope for the sum of two nonconvex functions:
  Further properties and nonmonotone linesearch algorithms,''
\newblock {\em SIAM Journal on Optimization}, vol. 28, no. 3, pp. 2274--2303,
  2018.

\bibitem{giselsson2018envelope}
Pontus Giselsson and Mattias F{\"a}lt,
\newblock ``Envelope functions: Unifications and further properties,''
\newblock {\em Journal of Optimization Theory and Applications}, vol. 178, no.
  3, pp. 673--698, 2018.

\bibitem{parikh_proximal_2014}
Neal Parikh and Stephen Boyd,
\newblock ``Proximal {Algorithms},''
\newblock {\em Foundations and Trends® in Optimization}, vol. 1, no. 3, pp.
  127--239, 2014.

\bibitem{combettes_proximal_2011}
Patrick~L. Combettes and Jean-Christophe Pesquet,
\newblock ``Proximal {Splitting} {Methods} in {Signal} {Processing},''
\newblock in {\em Fixed-{Point} {Algorithms} for {Inverse} {Problems} in
  {Science} and {Engineering}}, vol.~49, pp. 185--212. Springer New York, New
  York, NY, 2011.

\bibitem{popkov2005gradient}
A~Yu Popkov,
\newblock ``Gradient methods for nonstationary unconstrained optimization
  problems,''
\newblock {\em Automation and Remote Control}, vol. 66, no. 6, pp. 883--891,
  2005.

\bibitem{dontchev2013euler}
Asen~L Dontchev, MI~Krastanov, R~Tyrrell Rockafellar, and Vladimir~M Veliov,
\newblock ``An {Euler--Newton} continuation method for tracking solution
  trajectories of parametric variational inequalities,''
\newblock {\em SIAM Journal on Control and Optimization}, vol. 51, no. 3, pp.
  1823--1840, 2013.

\bibitem{dontchev2014implicit}
Asen~L Dontchev and R~Tyrrell Rockafellar,
\newblock {\em Implicit Functions and Solution Mappings: A View from
  Variational Analysis},
\newblock Springer, 2014.

\bibitem{meyer2000matrix}
Carl~D. Meyer,
\newblock {\em Matrix Analysis and Applied Linear Algebra}, vol.~71,
\newblock SIAM, 2000.

\end{thebibliography}
\end{document}